\documentstyle[12pt,a4,amssymb]{amsart}

\numberwithin{equation}{section}

\newtheorem{thm}{Theorem}[section]

\newtheorem{lem}{Lemma}[section]
\newtheorem{prop}{Proposition}[section]
\newtheorem{defi}{Definition}[section]
\newtheorem{rem}{Remark}[section]

\newcommand*{\C}{\mathbb{C}}
\newcommand*{\R}{\mathbb{R}}
\newcommand*{\Z}{\mathbb{Z}}
\newcommand*{\Q}{\mathbb{Q}}
\newcommand*{\A}{\mathbb{A}}

\begin{document}
\title[\null]
{
Two dimensional adelic analysis and \\ cuspidal automorphic representations of ${\rm GL}(2)$}
\author[\null]{Masatoshi Suzuki}
\date{\today \\ \quad This work was supported by JSPS}


\maketitle

\begin{abstract}
Two dimensional adelic objects were introduced by I. Fesenko
in his study of the Hasse zeta function associated to a regular model $\mathcal E$ of the elliptic curve $E$. 
The Hasse-Weil $L$-function $L(E,s)$ of $E$ appears in the denominator of the Hasse zeta function of $\mathcal E$. 
The two dimensional adelic analysis predicts that 
the integrand $h$ of the boundary term of the two dimensional zeta integral attached to $\mathcal E$ 
is mean-periodic. 
The mean-periodicity of $h$ implies 
the meromorphic continuation and the functional equation of $L(E,s)$. 
On the other hand, if $E$ is modular, several nice analytic properties of $L(E,s)$, 
in particular the analytic continuation and the functional equation,  
are obtained by the theory of the cuspical automorphic representation of $GL(2)$ 
over the ordinary ring of adele (one dimensional adelic object). 
In this article we try to relate the theory of two dimensional adelic object 
to the theory of cuspidal automorphic representation of $GL(2)$ over the one dimensional adelic object, 
under the assumption that $E$ is modular. 
Roughly speaking, they are dual each other.
\end{abstract}

\section{Introduction}
Let $X \to {\rm Spec}\,{\Z}$ be a scheme separated and of finite type.
The Hasse zeta function of $X$ is defined by the Euler product
\begin{equation*}
\zeta_X(s) = \prod_{ x \in X_0} \left(1-|\kappa(x)|^{-s} \right)^{-1},
\end{equation*}
where $X_0$ is the set of all closed points $x$ of $X$ with residue field $\kappa(x)$ of cardinality $|\kappa(x)|<\infty$.
For a number field $k$ with the ring of integers ${\mathcal O}_k$
the Hasse zeta function of the affine scheme ${\rm Spec} \, {\mathcal O}_k$
is the Dedekind zeta function $\zeta_k(s)=\prod_{{\frak p} \subset {\mathcal O}_k}(1-|{\mathcal O}_k/{\frak p}|^{-s})^{-1}$.
It is conjectured that $\zeta_X(s)$ has several nice analytic properties such as a meromorphic continuation and a functional equation.
However, the known result is very few when the dimension of $X$ is larger than one.

If the dimension of $X$ is one, the Hasse zeta function $\zeta_X(s)$ is essentially the Dedekind zeta funtion $\zeta_k(s)$.
Due to the celebrated work of Iwasawa and Tate,
the analytic properties of $\zeta_k(s)$ are obtained by the Fourier analysis on adele ${\A}_k$.
The completed Dedekind zeta function $\widehat{\zeta}_k(s)$ is defined by multiplying $\zeta_k(s)$ with a finite product of $\Gamma$-factors.
It has the integral representation
\begin{equation*}
\widehat{\zeta}_k(s)=\int_{{\A}_k^\times} f(x) |x|^s d \mu_{{\A}_k^\times}(x)=:\zeta_k(f,s),
\end{equation*}
where $f$ is an appropriate Schwartz--Bruhat function on ${\A}_k$ and $|~|$ is a module on the ideles ${\A}_k^\times$ of $k$.
On the other hand, one has
\begin{equation*}
\zeta_k(f,s)=\xi(f,s)+\xi(\widehat{f},1-s) + \omega(f,s)
\end{equation*}
on $\Re(s)>1$, where $\widehat{f}$ is the Fourier transform of $f$ on ${\A}_k$, $\xi(f,s)$ is an entire function given by an integral which converges absolutely for any $s \in \C$ and the boundary term
\begin{equation*}
\omega(f,s)=\int_{0}^1h_f(x)x^s \frac{dx}{x}
\end{equation*}
for some function $h_f$ on $(0,1)$.
The meromorphic continuation and the functional equation for $\widehat{\zeta}_k(s)$ are equivalent to the meromorphic
continuation and the functional equation for $\omega(f,s)$.
The properties of the function $h_f(x)$ are crucial in order to have a better understanding of $\omega(f,s)$.
Fourier analysis and analytic duality on $k \subset {\A}_k$ leads
\begin{equation*}
h_f(x)=-\mu\left({\A}_k^1 \slash k^\times \right) \left(f(0)-x^{-1}\widehat{f}(0)\right).
\end{equation*}
As a consequence, $\omega(f,s)$ is a rational function of $s$ invariant
with respect to $f\mapsto\widehat{f}$ and $s\mapsto(1-s)$.
Thus, $\widehat{\zeta}_k(s)$ admits a meromorphic continuation to $\C$
and satisfies a functional equation with respect to $s\mapsto(1-s)$.

Let $E$ an elliptic curve over $k$ and let ${\mathcal E} \to B={\rm Spec} {\mathcal O}_k$ be a regular model of $E$ over $k$.
Then the description of geometry of models in \cite[Thms 3.7, 4.35 in Ch. 9 and section 10.2.1 in Ch. 10]{Li} implies that
\begin{equation} \label{101}
\zeta_{\mathcal E}(s) = n_{\mathcal E}(s)\zeta_E(s) \quad \text{with} \quad
\zeta_E(s) = \frac{\zeta_k(s)\zeta_k(s-1)}{L(E,s)}
\end{equation}
on $\Re(s)>2$. Here $n_{\mathcal E}(s)$ is the product of zeta functions of affine lines
over finite extension $\kappa(b_j)$ of the residue fields $\kappa(b)$:
\begin{equation} \label{n_E}
n_{\mathcal E}(s) = \prod_{j=1}^J \left( 1 - |\kappa(b_j)|^{1-s} \right)^{-1}
\end{equation}
where $J$ is the number of singular fibres of ${\mathcal E} \to B$ (see \cite[section 7.3]{Fe3}).

The modularity conjecture for $E/k$ asserts that there exists a cuspidal automorphic representation $\pi_E$ of ${\rm GL}_2({\A}_k)$
such that
\[
L(E,s) = L(\pi_E,s-1/2).
\]
Then the general theory of $L$-function $L(\pi,s)$ of cuspidal automorphic representation $\pi$ of ${\rm GL}_2({\A}_k)$ leads to
an analytic continuation and a functional equation of $L(E,s)$ via $L(\pi,s)$.
The analytic properties of $L(\pi,s)$ are obtained by extending the Iwasawa-Tate theory
from the commutative group ${\rm GL}_1({\A}_k)$ to the noncommutative group ${\rm GL}_2({\A}_k)$.
In this story, the theory of noncommutative group ${\rm GL}_2({\A}_k)$ relates to $\zeta_{\mathcal E}(s)$
via the modularity conjecture and the \emph{$L$-function $L(E,s)$ of $E$}.

In contrast with the above story, I. Fesenko proposed another way to study $\zeta_{\mathcal E}(s)$ in \cite{Fe1,Fe3,Fe2}
by using a commutative group associated to two dimensional adeles. The ordinary ring of adeles ${\A}_k$ is regarded as an one dimensional object
in the sense that it is associated to the one dimensional scheme ${\rm Spec} \, {\mathcal O}_k$.
He introduced the two dimensional adelic space ${\bf A}_{\mathcal E}$ associated to the two dimensional scheme $\mathcal E$
and established a theory of translation invariant measure and integrals on its subring ${\A}_{{\mathcal E},S} \prec {\bf A}_{\mathcal E}$,
where $S$ is a set of fibers consisting of finitely many horizontal curves of ${\mathcal E} \to B$
and all its vertical fibers. Using a measure theory on two dimensional adelic space,
he defined the zeta integral
\[
\zeta_{{\mathcal E},S}(f,s) = \int_{T_{{\mathcal E},S}} f(t) |t|^s d\mu(t),
\]
where $f$ is a generalized Schwartz-Bruhat function on ${\A}_{{\mathcal E},S} \times {\A}_{{\mathcal E},S}$,
$T_{{\mathcal E},S}$ is certain subgroup of ${\A}_{{\mathcal E},S}^\times \times {\A}_{{\mathcal E},S}^\times$,
$|\,\,|$ is a module function on $T_{{\mathcal E},S}$ and $d\mu$ is a measure on $T_{{\mathcal E},S}$ (see \cite[section 5]{Fe3}).
The zeta integral $\zeta_{{\mathcal E},S}(f,s)$ converges absolutely for $\Re(s)>2$ .
If the test function $f_0$ is well-chosen, the zeta integral $\zeta_{{\mathcal E},S}(f_0,s)$ equals
\[
\zeta_{{\mathcal E},S}(f_0,s)
= \prod_{\rm finite} \widehat{\zeta}_{k_i}(s/2)^2 \cdot c_{\mathcal E}^{1-s} \cdot \zeta_{\mathcal E}(s)^2
\]
where $k_i$ is an extension of $k$ determined by each horizontal fiber in $S$
and $c_{\mathcal E}$ is a positive real number determined by $\mathcal E$.
On the other hand, similar to the Iwasawa-Tate theory, the zeta integral $\zeta_{{\mathcal E},S}(f,s)$
is decomposed as
\[
\zeta_{{\mathcal E},S}(f,s) = \xi(f,s)+\xi(\widehat{f},2-s) + \omega(f,s)
\]
on $\Re(s)>2$, where $\xi(f,s)$ is an entire function and
$\widehat{f}$ is the Fourier transform of $f$ on ${\A}_{{\mathcal E},S} \times {\A}_{{\mathcal E},S}$.
Hence the meromorphic continuation of $\omega(f_0,s)$
implies the meromorphic continuation of the Hasse zeta function $\zeta_{\mathcal E}(s)$.
If we can prove the meromorphic continuation of $\omega(f,s)$
using analysis and duality on two dimensional adelic space ${\A}_{{\mathcal E},S}$,
it leads the meromorphic continuation of the $L$-function $L(E,s)$, without proving the modularity property!

One possible approach for the meromorphic continuation of $\omega(f,s)$
is proposed via the theory of mean-periodic functions (\cite[section 7]{Fe3}, see also \cite{IGM}).
For the general theory of mean-periodic functions, see Kahane \cite{Kah}, Schwartz \cite{MR0023948} or a reference of \cite{IGM}.
Similar to the Iwasawa-Tate theory, we have the boundary term
\[
\omega(f,s) = \int_0^1 h_f(x) \cdot x^{s} \frac{dx}{x} = \int_0^\infty h_f(e^{-t}) \cdot  e^{-st}dt
\]
for some function $h_f$ on $(0,1)$.
So the boundary term is the Laplace transform of $h_f(e^{-t})$.

Let $\mathfrak X$ be a locally convex separated topological $\C$-vector space
consisting of complex valued functions on ${\R}_+^\times = (0,\infty)$.
It has a natural representation $\tau$ of ${\R}_+^\times$ as
$(\tau_a F)(x)=F(x/a)$ for every $F \in \mathfrak X$.
For $F \in \mathfrak X$ we denote by  ${\mathcal T}(F)$ be the closure of
$\{\tau_a F \,|\, a \in {\R}_+^\times \}$ with respect to the topology of $\mathfrak X$.
A function $F \in \mathfrak X$ is called mean-periodic if ${\mathcal T}(F) \not= \mathfrak X$.
Using the representation $\tau$ the convolution $F \ast \varphi$ for $F \in {\mathfrak X}$ and $\varphi \in {\mathfrak X}^\ast$
is defined by
\[
(F \ast \varphi)(x) = \langle \tau_x \check{F} ,\varphi \rangle
\]
where $\check{F}(x)=F(x^{-1})$.
The mean-periodicity ${\mathcal T}(F) \not= \mathfrak X$ of $F$
is equivalent that the space of annihilators $\varphi \in {\mathfrak X}^\ast$ of ${\mathcal T}(F)$
concerning the convolution is nontrivial;
\[
{\mathcal T}(F)^{\bot} := \{ \varphi \in X^\ast \,|\,  G \ast \varphi = 0,~
\forall G \in {\mathcal T}(F) \} \not= \{0\}.
\]
As a consequence of the general theory of mean-periodic function,
if $F$ is mean-periodic, the Laplace transform of $F(e^{-t})$ (the Mellin transform of $F(x)$)
is continued meromorphically to the whole complex plane.

Now we suppose that $h_{f_0} \in \mathfrak X$.
Then the conjectural mean-periodicity of $h_{f_0}$ implies the meromorphic continuation of the Hasse zeta function $\zeta_{\mathcal E}(s)$.
Hence it is important to understand the space of annihilators ${\mathcal T}(h_{f_0})^{\bot}$.
\smallskip

In this paper, in the case $k={\Q}$, we describe the space of annihilators ${\mathcal T}(h_{f_0})^{\bot}$
by using the cuspidal automorphic representation of ${\rm GL}_2({\A}_{\Q})$ whose existence follows from the modularity of $E/\Q$
(see Theorem \ref{thm_01}, Theorem \ref{thm_02} for more detail).
Such description of ${\mathcal T}(h_{f_0})^{\bot}$ suggests some duality between the commutative theory of two dimensional adeles ${\bf A}_{\mathcal E}$, ${\A}_{{\mathcal E},S}$ and the noncommutative theory ${\rm GL}_2({\A}_{\Q})$ of one dimensional adele ${\A}_{\Q}$.
\bigskip

In section 2 we include several definitions, notations and already known properties.
In section 3 we state the results, and we prove them in section 4.

\section{Preliminaries}

Let $S(\R)$ be the Schwartz space on $\R$ which consists of smooth functions on $\R$ satisfying
\begin{equation*}
\Vert f \Vert_{m,n} = \sup_{x \in \R} \vert x^m f^{(n)}(x) \vert <\infty
\end{equation*}
for all nonnegative integer $m$ and $n$.
It is a Fr\'echet space over the complex numbers with the topology induced from the family of seminorms $\Vert ~\Vert_{m,n}$.
Let us define the Schwartz space $S(\R_+^\times)$ on ${\R}_+^\times$ and its topology via the homeomorphism
\begin{equation*}
S({\R}) \rightarrow  S({\R}_+^\times); ~f(t)  \mapsto  f(-\log{x}),
\end{equation*}
where $t = -\log x$.
The \emph{strong Schwartz space} $\mathbf{S}({\R}_+^\times)$ (\cite{Mey1}) is defined by
\begin{equation}\label{eq_strong}
{\mathbf S}({\R}_+^\times):=\bigcap_{\beta\in\R}
\left\{f:{\R}_+^\times\to{\C},\left[x \mapsto x^{-\beta}f(x)\right] \in S(\R_+^\times)\right\}.
\end{equation}
One of the family of seminorms on $\mathbf{S}({\R}_+^\times)$ defining its topology is given by
\begin{equation} \label{901}
\Vert f \Vert_{m,n}=\sup_{x \in {\R}_+^\times} \vert x^m f^{(n)}(x) \vert
\end{equation}
for integer $m$ and nonnegative integer $n$.
The strong Schwartz space $\mathbf{S}({\R}_+^\times)$ is a Fr\'echet space over the complex numbers
where the family of seminorms defining its topology is given in \eqref{901}.
This space is closed under the multiplication by a complex number and the pointwise addition and multiplication (\cite{Mey1}).
Let $\mathbf{S}(\R_+^\times)^\ast$ be the dual space of $\mathbf{S}(\R_+^\times)$ with the weak $\ast$-topology.
The pairing between $\mathbf{S}(\R_+^\times)$ and $\mathbf{S}(\R_+^\times)^\ast$ is denoted by $\langle ~, ~\rangle$, namely
$\langle f,\varphi\rangle=\varphi(f)$ for  $f\in\mathbf{S}(\R_+^\times)$ and  $\varphi\in\mathbf{S}(\R_+^\times)^\ast$.
The (multiplicative) representation $\tau$ of $\R_+^\times$ on $\mathbf{S}({\R_+^\times})$ is defined by
\begin{equation*}
\tau_x f(y):=f(y/x), \quad \forall x\in{\R}_+^\times
\end{equation*}
and the (multiplicative) convolution $f \ast \varphi$ of
$f\in\mathbf{S}({\R_+^\times})$ and $\varphi\in\mathbf{S}(\R_+^\times)^\ast$ by
\begin{equation*}
(f \ast \varphi)(x)= \langle \tau_x \check{f},\varphi\rangle,
\quad \forall x \in {\R}_+^\times
\end{equation*}
where $\check{f}(x):=f(x^{-1})$.
The dual representation $\tau^{\ast}$ on ${\mathbf S}(\R_+^\times)^\ast$ is defined by
\begin{equation*}
\langle f,\tau_x^{\ast}\varphi\rangle:=\langle\tau_x f,\varphi\rangle.
\end{equation*}
If $V$ is a $\C$-vector space then the bidual space $V^{\ast\ast}$
(the dual space of $V^\ast$ with respect to the weak $\ast$-topology on $V^\ast$)
is identified with $V$ in the following way.
For a continuous linear functional $F$ on $V^\ast$ with respect to its weak $\ast$-topology,
there exists $v \in V$ such that $F(v^\ast)=v^\ast(v)$ for every $v^\ast \in V^\ast$.
Therefore, we do not distinguish the pairing on $V^{\ast\ast}\times V^{\ast}$ from the pairing on $V\times V^\ast$.
\begin{defi}
Let $\mathfrak{X} = \mathbf{S}({\R_+^\times})^\ast$.
An element $x\in\mathfrak{X}$ is said to be $\mathfrak{X}$-mean-periodic
if there exists a non-trivial element $x^\ast$ in $\mathfrak{X}^\ast$ satisfying $x \ast  x^\ast=0$.
\end{defi}
For $x \in \mathfrak{X} = \mathbf{S}({\R_+^\times})^\ast$, we denote by ${\mathcal T}(x)$ the closure of the $\C$-vector space spanned by
$\{\tau_g^\ast(x), g \in {\R}_+^\times \}$.
The Hahn-Banach theorem leads to another definition of $\mathfrak{X}$-mean-periodic functions.
\begin{prop}\label{proposition_td_mp}%
An element $x\in \mathfrak{X}=\mathbf{S}({\R_+^\times})^\ast$ is $\mathfrak{X}$-mean-periodic
if and only if ${\mathcal T}(x)\not=\mathfrak{X}$.
\end{prop}

Let $L_{\textrm{loc},{\rm poly}}^1({\R_+^\times})$
be the space of locally integrable functions on ${\R_+^\times}$ satisfying
\begin{equation*}
h(x)=\begin{cases}
O(x^{a}) & \text{as $x\to+\infty$}, \\
O(x^{-a}) & \text{as $x\to0^+$}
\end{cases}
\end{equation*}
for some real number $a\geq 0$. Each $h\in L_{\textrm{loc},{\rm poly}}^1({\R_+^\times})$ gives rise to a distribution $\varphi_h \in \mathbf{S}({\R_+^\times})^\ast$ defined by
\begin{equation*}
\langle f , \varphi_h \rangle=\int_{0}^{+\infty}f(x)h(x)\frac{dx}{x},
\quad \forall f \in {\mathbf S}({\R}_+^\times).
\end{equation*}
If there is no confusion, we denote $\varphi_h$ by $h$ itself and use the notations $\langle f,h\rangle=\langle f,\varphi_h\rangle$ and $h(x)\in \mathbf{S}({\R}_+^\times)^\ast$.
Then
\begin{equation*}
x^\lambda\log^k{(x)}\in L_{\textrm{loc},{\rm poly}}^1({\R_+^\times})\subset \mathbf{S}({\R}_+^\times)^\ast
\end{equation*}
for all $k\in{\Z}_{ \ge 0}$ and $\lambda\in{\C}$.
Moreover, if $h\in L_{\textrm{loc},{\rm poly}}^1({\R_+^\times})$
then the convolution $f\ast\varphi_h$ coincides with the ordinary convolution on functions on ${\R}_+^\times$ namely
\[
(f \ast h)(x)
= \langle \tau_x \check{f},f\rangle
= \int_{0}^{+\infty}f(x/y)h(y)\frac{dy}{y}
= \int_{0}^{+\infty}f(y)h(x/y)\frac{dy}{y}.
\]
For a $h \in L_{\textrm{loc},{\rm poly}}^1({\R_+^\times})$ define $h^+$ and $h^-$ by
\begin{equation*}
h^+(x):=\begin{cases}
0 & \text{if $x\geq 1$}, \\
h(x) & \text{otherwise}
\end{cases}
\quad
h^-(x):=\begin{cases}
h(x) & \text{if $x\geq 1$}, \\
0 & \text{otherwise}.
\end{cases}
\end{equation*}
Clearly, $h^{\pm}\in L_{\textrm{loc},{\rm poly}}^1({\R_+^\times})$ for all  $h\in L_{\textrm{loc},{\rm poly}}^1({\R_+^\times})$.
\begin{lem}\label{lem_01}%
Let $h \in L_{\text{\emph{loc}},{\rm poly}}^1({{\R}_+^\times})$.
If $f \ast h=0$ for some non-trivial $f\in{\mathbf S}({{\R}_+^\times})$ then the Mellin transforms
\begin{equation*}
{\mathsf M}(f \ast h^{\pm})(s)=\int_{0}^{+\infty}(f \ast h^{\pm})(x)x^s\frac{dx}{x}
\end{equation*}
are entire functions on ${\C}$.
\end{lem}
\begin{defi}\label{def_MC}
Let $h\in L_{\text{\emph{loc}},{\rm poly}}^1({\R_+^\times})$.
If $f \ast h=0$ for some non-trivial $f\in\mathbf{S}({\R_+^\times})$
then the \emph{Mellin--Carleman transform} ${\mathsf MC}(h)(s)$ of $h(x)$ is defined by
\begin{equation*}
{\mathsf MC}(h)(s)
:=\frac{{\mathsf M}(f\ast h^+)(s)}{{\mathsf M}(f)(s)}
 =-\frac{{\mathsf M}(f\ast h^-)(s)}{{\mathsf M}(f)(s)}.
\end{equation*}
\end{defi}
The Mellin--Carleman transform ${\mathsf MC}(h)$
does not depend on the particular choice of non-trivial $f$ satisfying $f \ast h=0$.
By Lemma \ref{lem_01} we have
\begin{prop}
Let $h\in L_{\text{\emph{loc}},{\rm poly}}^1({\R_+^\times}) \subset {\mathbf S}({\R_+^\times})^\ast$.
If $f \ast h=0$ for some non-trivial $f\in{\mathbf S}({\R_+^\times})$, in other words,
$h$ is ${\mathbf S}({\R_+^\times})^\ast$-mean-periodic,
then the Mellin--Carleman transform ${\mathsf MC}(h)(s)$ of $h(x)$
is a meromorphic function on $\C$.
\end{prop}
The Mellin--Carleman transform ${\mathsf MC}(h)(s)$ of $h(x)$ is \emph{not} a generalization of the Mellin transform of $h$
but is a generalization of the following integral, half Mellin transform,
\[
\int_{0}^{1}h(x)x^{s}\frac{dx}{x}.
\]
See also section 2 of \cite{IGM} for more detail.
\smallskip

Let $E$ be an elliptic curve over $\Q$ with conductor $q_E$.
Then the completed $L$-function $\Lambda(E,s)$ is defined by
\[
\Lambda(E,s) := q_E^{s/2}(2\pi)^{-s}\Gamma(s)L(E,s).
\]
It is conjectured that $\Lambda(E,s)$ is continued to an entire function and satisfies
the functional equation
$\Lambda(E,s) = \omega_E \Lambda(E,2-s)$ for some sign $\omega_E \in \{\pm 1\}$.
By \eqref{101}, the meromorphic continuation and the functional equation of $\Lambda(E,s)$
implies the meromorphic continuation and the functional equation of $\zeta_{\mathcal E}(s)$.
Moreover such nice analytic properties of $\Lambda(E,s)$ lead to mean-periodicity of the $\omega(f_0,s)$.

\begin{thm} \label{prop_01}
Let $E$ be an elliptic curve over $\Q$ and let ${\mathcal E} \to {\rm Spec}\,{\Z}$ be its regular model.
Assume that $\Lambda(E,s)$ is continued meromorphically to ${\C}$ with a finite poles
and satisfies the functional equation
\[
\Lambda(E,s)^2 = \Lambda(E,2-s)^2.
\]
Then the function
\[
h_{\mathcal E}(x):=f_{\mathcal E}(x) - x^{-1} f_{\mathcal E}(x^{-1})
\]
with
\[
f_{\mathcal E}(x)
= \frac{1}{2\pi i} \int_{(c)}
\Lambda(s/2+1/4)^2 c_{\mathcal E}^{-s-1/2} \zeta_{\mathcal E}(s+1/2)^2 x^{-s} ds \quad (c>1)
\]
belongs to ${\mathbf S}({\R_+^\times})^\ast$,
where $c_{\mathcal E}$ is a positive real constant determined by the singular fiber of $\mathcal E$ {\rm (\cite[section 5]{Fe3})}.
Moreover $h_{\mathcal E}$ is ${\mathbf S}({\R_+^\times})^\ast$-mean-periodic and has the expansion
\[
\aligned
h_{\mathcal E}(x)
& = \lim_{T\to\infty} \sum_{\Im(\lambda)\leq T} \sum_{m=1}^{m_\lambda} C_m(\lambda) \frac{(-1)^{m-1}}{(m-1)!} x^{-\lambda} (\log x)^{m-1}
\endaligned
\]
where $\lambda$ are poles of $\Lambda(s/2+1/4)^2 c_{\mathcal E}^{-s-1/2} \zeta_{\mathcal E}(s+1/2)^2$ of multiplicity $m_\lambda$,
$C_m(\lambda)$ are constants determined by the principal part at $s=\lambda$;
\[
\Lambda(s/2+1/4)^2 c_{\mathcal E}^{-s-1/2} \zeta_{\mathcal E}(s+1/2)^2
= \sum_{m=1}^{m_\lambda} \frac{C_m(\lambda)}{(s-\lambda)^m} + O(1) \quad \text{when} \quad s \to \lambda,
\]
and the sum over $\lambda$ is converges uniformly on every compact subset of ${\R}_+^\times$.
\end{thm}
\begin{pf}
See section 5 of \cite[section 5]{IGM}.
\end{pf}
So the mean-periodicity of $h_{\mathcal E}(x)$ and the meromorphic continuation of $\Lambda(E,s)^2$
are equivalent to each other in the first approximation.
\begin{rem} Let $S$ be the set of fibres of ${\mathcal E} \to {\rm Spec}\,{\Z}$
consisting of one horizontal curve which is the image of the zero section of ${\mathcal E} \to {\rm Spec}\,{\Z}$
and all vertical fibres of ${\mathcal E} \to {\rm Spec}\,{\Z}$.
Then we have
\[
\aligned
\zeta_{{\mathcal E},S}(f_0,s)
& = \Lambda(s/2)^2 c_{\mathcal E}^{-s} \zeta_{\mathcal E}(s)^2 \\
& = \int_{1}^{\infty} x^{-1/2} f_{\mathcal E}(x) x^{s} \frac{dx}{x}
+ \int_{1}^{\infty} x^{-1/2} f_{\mathcal E}(x) x^{2-s} \frac{dx}{x}
+ \int_{0}^{1} x^{-1/2} h_{\mathcal E}(x) x^{s} \frac{dx}{x}.
\endaligned
\]
Hence the function $h_{f_0}(x)$ in the introduction is $x^{-1/2} h_{\mathcal E}(x)$.
\end{rem}
\begin{rem}
We hope to prove the mean-periodicity of $h_{\mathcal E}(x)$ without using the meromorphic continuation of $\Lambda(E,s)$.
\end{rem}

\section{Statement of Results}
Throughout this section we denote by $\A$ the adele ${\A}_{\Q}$ of ${\Q}$.
At first we settle the following basic assumption.
\smallskip

\noindent
{\bf Basic assumption.} {\it Suppose that $E/{\Bbb Q}$ is modular. We denote by $(\pi,V_\pi)$
the corresponding cuspical automorphic representation
in $L^2(GL_2({\Q}) \backslash GL_2({\A}),{\mathbf 1})$,
where ${\mathbf 1}$ is the trivial central character. }
\smallskip

Of course the modularity of $E/{\Bbb Q}$ is now a theorem by the famous work of Wiles et. al.
However it is not proved for a general number field $k$.
We emphasize this assumption for the future study of this direction.

\subsection{Construction on the positive real line}
In this part we construct the space of anihilators ${\cal T}(h_{\mathcal E})^{\bot}$
of $h_{\mathcal E}$ associated to $\zeta_{\mathcal E}(s)$ as in Theorem \ref{prop_01}
by using ${\rm GL}_2({\Bbb A})$-theory of Soul\'e ~\cite{So}
which is an extension of the original theory of Connes \cite{Co}.
\smallskip

Let $M={\rm Mat}_2$ and $G={\rm GL}_2$.
Let $|~| : G_{\Bbb A} \to {\Bbb R}_+^\times$ be
the module map given by $|g|=|\det g|_{\A}$.
Let $f_\pi$ be an admissible matrix coefficient
of the cuspidal automorphic representation $(\pi,V_\pi)$ on $L^2(G_{\Bbb Q} \backslash G_{\Bbb A},{\mathbf 1})$,
and let $\phi$ be a Schwartz-Bruhat function on $M_{\Bbb A}$.
For $x \in {\R}_+^\times$, we set $G_x = \{ g \in G_{\Bbb A} \,|\, |g|=x\}$.
Define a complex valued function ${\frak E}(\phi,f_\pi)$ on ${\R}_+^\times$ by
\begin{equation} \label{801}
{\frak E}(\phi,f_\pi)(x) = \int_{G_x} \phi(g) f_\pi(g) dg \quad (x \in {\R}_+^\times).
\end{equation}
Then
\begin{enumerate}
\item[i)] the integral \eqref{801} converges absolutely,
\item[ii)] for any integer $N>0$,
there exists a positive constant $C$ such that
\begin{equation} \label{802}
\vert {\frak E}(\phi,f_\pi)(x) \vert \leq C x^{-N},
\end{equation}
for all $x \in {\R}_+^\times$,
\item[iii)] we have the functional equation
\begin{equation} \label{803}
{\frak E}(\phi,f_\pi)(x) = x^{-2} {\frak E}(\hat{\phi}, \check{f}_\pi)(x^{-1})
\end{equation}
where $\hat{\phi}$ is the Fourier transform of $\phi$ and $\check{f}_\pi(g)=f_\pi(g^{-1})$.
\end{enumerate}
Let
\[
S(\pi) = \{(\phi,f_\pi) \,|\, \phi \in S(M({\Bbb A})),~\text{$f_\pi$: admissible coefficient of $\pi$} \}.
\]
Then \eqref{802} and \eqref{803} show that ${\frak E}$ is a map from $S(\pi)$ into ${\bf S}({\Bbb R}_+^\times)$:
\[
{\frak E} : S(\pi) \to {\bf S}({\Bbb R}_+^\times); ~
(\phi,f_\pi) \mapsto {\frak E}(\phi,f_\pi).
\]
We denote by ${\cal V}_\pi  \subset {\bf S}({\Bbb R}_+^\times)$ the image of ${\frak E}$.
Using the function
\begin{equation} \label{804}
w_0(x) = \frac{1}{2\pi i} \int_{(c)} \Gamma(s/4)^2
\cdot \frac{(c_{\mathcal E}/q_E)^{s}}{n_{\mathcal E}(s)^2}
\cdot s^4 (s-2)^4 \cdot (s-1)^2 \cdot x^{-s}ds,
\end{equation}
we define the space ${\cal W}_\pi$ by
\begin{equation} \label{805}
{\cal W}_\pi = w_0 \ast {\cal V}_\pi \ast {\cal V}_\pi =
{\rm span}_{\Bbb C}\{ w_0 \ast v_1 \ast v_2 \,|\, v_i \in {\cal V}_\pi \}.
\end{equation}
Then the space ${\cal W}_\pi$ is a subspace of ${\bf S}({\Bbb R}_+^\times)$,
since $w_0 \in {\bf S}({\Bbb R}_+^\times)$ and ${\bf S}({\Bbb R}_+^\times)$ is closed under the multiplicative convolution.
For $h \in {\bf S}({\R}_+^\times)$, we define
\[
\aligned
{\mathcal T}(h)^\perp
& = \{\, g \in {\bf S}({\R}_+^\times) \,|~g \ast \tau =0, ~\forall \tau \in {\mathcal T}(h) \}
\endaligned
\]
and
\[
{\cal W}_\pi^\perp
= \{\, \varphi \in {\bf S}({\Bbb R}_+)^\ast \,|\, w \ast \varphi =0, ~\forall w \in {\cal W}_\pi \,\}.
\]
\begin{thm} \label{thm_01}
Let $E$ be an elliptic curve over $\Q$ and let ${\mathcal E} \to {\rm Spec}\, {\Z}$ be its regular model.
Let $h_{\mathcal E}$ be the function on ${\R}_+^\times$ associated to the Hasse zeta function $\zeta_{\mathcal E}(s)^2$
as in Theorem \ref{prop_01}.
Then we have
\[
{\cal T}(h_{\mathcal E})^\perp  \supset {\cal W}_\pi
\quad \text{and} \quad
{\cal T}(h_{\mathcal E})  \subset {\cal W}_\pi^\perp.
\]
Hence ${\cal W}_{\pi} \not=\{0\}$ means ${\bf S}({\Bbb R}_+)^\ast$-mean-periodicity of $h_{\mathcal E}(x)$.
Further the equality
\[
{\cal T}(h_{\mathcal E})^\perp = {\cal W}_\pi
\quad \text{or} \quad
{\cal T}(h_{\mathcal E}) = {\cal W}_\pi^\perp
\]
implies the non-existence of cancelations of zeros
between \[
\text{$(s-1) \widehat{\zeta}(s/2) \widehat{\zeta}(s) \widehat{\zeta}(s-1)$
\quad and \quad $n_{\mathcal E}(s)^{-1}\Lambda_E(s)$.}
\]
\end{thm}

\subsection{Adelic construction}
In this part, we consider the adelic version of the previous one
according to Deitmar~\cite{De}.
Let $S(M_{\Bbb A})_0$ be the space of all $\phi \in S(M_{\Bbb A})$
such that $\phi$ and $\hat{\phi}$ send
$\{g \in M_{\Bbb A} \,|\, \det(g)=0\} = M_{\Bbb A} - G_{\Bbb A}$ to zero.
For $\phi \in S(M_{\Bbb A})_0$
we define functions ${\frak E}(\phi)$ and $\hat{\frak E}(\phi)$ on $G_{\Bbb A}$ by
\begin{equation} \label{807}
\aligned
{\frak E}(\phi)(g) & = \sum_{\gamma \in M_{\Bbb Q}} \phi(\gamma g) = \sum_{\gamma \in G_{\Bbb Q}} \phi(\gamma g), \\
\hat{\frak E}(\phi)(g) & = \sum_{\gamma \in M_{\Bbb Q}} \phi( g \gamma) = \sum_{\gamma \in G_{\Bbb Q}} \phi(g \gamma).
\endaligned
\end{equation}
Then for any $\phi \in S(M_{\Bbb A})_0$, we have
\begin{enumerate}
\item[i)] the sums ${\frak E}(\phi)$ and $\hat{\frak E}(\phi)$
converge locally uniformly in $g$ with all derivatives,
\item[ii)] for any $N>0$ there exists $C>0$ such that
\begin{equation} \label{808}
|{\frak E}(\phi)(g)|,\,|\hat{\frak E}(\phi)(g)| \leq C {\rm min}\bigl( |g|, |g|^{-1} \bigr)^N,
\end{equation}
\item[iii)] for $g \in G_{\Bbb A}$ we have the functional equation
\begin{equation} \label{809}
{\frak E}(\phi)(g) = |g|^{-2} \hat{\frak E}(\hat{\phi})(g^{-1}).
\end{equation}
\end{enumerate}
Hence ${\frak E}(\phi)$ belongs to the strong Schwartz space
\[
{\bf S}(G_{\Q} \backslash G_{\A} ) =
\bigcap_{\beta \in {\R}} |~|^\beta S( G_{\Q} \backslash G_{\A} ).
\]
Let $G_{\A}^1$ be the kernel of the module map $g \mapsto |g|$.
Fix a splitting ${\ss} : {\R}_+^\times \to G_{\A}$ of the exact sequence $1 \to G_{\A}^1 \to G_{\A} \to 1$
such that $({\rm id},{\ss}):G_{\A}^1 \times {\R}_+^\times \to G_{\A}$ is an isomorphism.
We denote by ${\mathsf R}$ the image of splitting $\ss$.
Let $\varphi_\pi \in V_\pi \subset L^2(G_{\Q} \backslash G_{\A}^1) \simeq L^2({\mathsf R} G_{\Q} \backslash G_{\Bbb A})$
be a vector $\varphi_\pi = \otimes_v \varphi_{\pi,v}$
such that $\varphi_{\pi,v}$ is a normalized class one vector
for almost all places. Further we assume that
$\varphi_\pi$ is smooth and $\varphi_\pi(1) \not=0$.

We define
\[
{\mathcal W}_\pi =
{\rm span}_{\C}
\bigl\{
(w_0 \circ \ss) \ast
({\frak E}(\phi_1) \cdot \varphi_\pi) \ast
({\frak E}(\phi_2) \cdot \varphi_\pi) \,\bigl|\, \phi_i \in S(M_{\A})
\bigr\} \subset {\bf S}( G_{\Q} \backslash G_{\A}),
\]
where $w_0$ is the function in \eqref{804},
$({\frak E}(\phi)\cdot \varphi_\pi)(x)={\frak E}(\phi)(x) \cdot \varphi_\pi(x)$
and $\ast$ is the convolution on $G_{\Bbb Q} \backslash G_{\Bbb A}$
via the right regular representation $R$, and
\[
{\cal W}_\pi^{\perp}
= \{ \eta \in {\bf S}( G_{\Q} \backslash G_{\A})^\ast
\,|\, w \ast \eta  \equiv 0, ~\forall w \in {\mathcal W}_\pi  \}.
\]
For $\eta \in {\bf S}( G_{\Bbb Q} \backslash G_{\Bbb A})^\ast$ we define
\[
{\cal T}(\eta)= {\rm span}_{\Bbb C}\{ R^\ast(g)\eta \,|\, g \in G_{\A} \},
\]
where $R^\ast$ is the transpose of the right regular representation of $G_{\A}$
on ${\bf S}( G_{\Q} \backslash G_{\A})$
with respect to the pairing $\langle ~,~ \rangle$ of ${\bf S}( G_{\Q} \backslash G_{\A})$ and ${\bf S}( G_{\Q} \backslash G_{\A})^\ast$,
\begin{thm} \label{thm_02}
Let $h_{\mathcal E}$ be the function on ${\R}_+^\times$ associated to the Hasse zeta function $\zeta_{\mathcal E}(s)^2$
as in Theorem \ref{prop_01}.
Under the above notations, we have
\[
{\cal T}(h_{\mathcal E} \circ {\ss}) \subset {\cal W}_\pi^\perp.
\]
The equality
\[
{\cal T}(h_{\mathcal E} \circ {\ss}) = {\cal W}_\pi^\perp
\]
implies the non-existence of cancelation of the zeros
between
\[
(s-1) \widehat{\zeta}(s/2) \widehat{\zeta}(s) \widehat{\zeta}(s-1)
\quad \text{and} \quad n_{\mathcal E}(s)^{-1}\Lambda_E(s). \]
\end{thm}

\section{Proof of Results}

\subsection{Proof of Theorem \ref{thm_01}}

First we prove the implication ${\cal W}_\pi \subset {\cal T}(h_{\mathcal E})^\perp$.
It suffices to prove that $w \ast h_{\mathcal E}\equiv 0$ for any $w \in {\cal W}_\pi.$
By Theorem \ref{prop_01} the function $h_{\mathcal E}$
is a series consisting of functions $f_{\lambda,k}(x)=x^{-\lambda}(\log x)^k$.
For $w(x) \in {\mathcal W}_\pi$,
\begin{equation*}
\aligned
w \ast f_{\lambda,k}(x)
& = \int_{0}^{\infty} w(y) f_{\lambda,k}(x/y) \frac{dx}{x} \\
& = \sum_{j=1}^k (-1)^k \binom kj x^{-\lambda}(\log x)^{k-j}
\int_{0}^{\infty} w(y) y^{\lambda} (\log y)^j \frac{dy}{y}.
\endaligned
\end{equation*}
Here
\begin{equation*}
\aligned
\int_{0}^{\infty} w(y) y^{\lambda} (\log y)^j \frac{dy}{y}
= \frac{d^j}{d\lambda^j} \int_{0}^{\infty} w(y) y^{\lambda} \frac{dy}{y}.
\endaligned
\end{equation*}
By definition of ${\mathcal W}_\pi$,
\begin{equation*}
\aligned
\int_{0}^{\infty} w(y) y^{\lambda} \frac{dy}{y}
& = \int_{0}^{\infty}
\left(
w_0 \ast {\frak E}(\phi_1,f_\pi) \ast {\frak E}(\phi_2,f_\pi^\prime)
\right)(y)
 y^{\lambda} \frac{dy}{y} \\
& =
\int_{0}^{\infty} w_0(y) y^{\lambda} \frac{dy}{y} \cdot
\int_{0}^{\infty} {\frak E}(\phi_1,f_\pi)(y) y^{\lambda} \frac{dy}{y} \cdot
\int_{0}^{\infty} {\frak E}(\phi_2,f_\pi^\prime)(y) y^{\lambda} \frac{dy}{y}
\endaligned
\end{equation*}
for some $(\phi_1,f_\pi),\, (\phi_2,f_\pi^\prime) \in S(\pi)$.
From the construction of ${\mathcal V}_\pi$, we have
\begin{equation*}
\int_{0}^{\infty} {\frak E}(\phi,f_\pi)(y) y^{\lambda} \frac{dy}{y}
= F_{\phi,f_\pi}(\lambda) L(\pi,\lambda-1/2) =  F_{\phi,f_\pi}(\lambda) L(E,\lambda)
\end{equation*}
where $F_{\phi,f_\pi}(\lambda)$ is an entire function determined by $(\phi,f_\pi)$ (see \cite[Theorem 13.8]{GoJa} and \cite[section 2.5]{So}).
The second equality is a consequence of modularity.
Therefore
\begin{equation} \label{401}
\aligned
\int_{0}^{\infty} w(y) y^{\lambda} \frac{dy}{y}
& = \Gamma(\lambda/4)^2 \lambda^4(\lambda-2)^4 (\lambda-1)^2 (c_{\mathcal E}/q_E)^{\lambda} \\
& \quad \times n_{\mathcal E}(\lambda)^{-2} L(E,\lambda)^2 F_{\phi_2,f_\pi}(\lambda) F_{\phi_2,f_\pi^\prime}(\lambda),
\endaligned
\end{equation}
and $w \ast h_{\mathcal E}$ is a series consisting of \eqref{401} and its $j$-th derivative with $j \leq m_\lambda$.
Because $E/{\Q}$ is modular, $\Lambda(E,s)$ is an entire function.
Therefore the complex numbers $\lambda$ appearing in the expansion
of $h_{\mathcal E}(x)$ is one of the followings:
\begin{enumerate}
\item $\lambda=0$ or $2$ and $m_\lambda = 4$,
\item $\lambda\not=1$ is a zero of $\Lambda(E,s)$ with $n_{\mathcal E}(\lambda)^{-1}\not=0$,
and $0 \le m_\lambda \le$ the multiplicity of zero of $\Lambda(E,s)^2$ at $s=\lambda$,
\item $\lambda\not=1$ is a common zero of $\Lambda(E,s)$ and $n_{\mathcal E}(s)^{-1}$,
and $-2 \le m_\lambda-2 \le$ the multiplicity of zero of $\Lambda(E,s)^2$ at $s=\lambda$,
\item $\lambda\not=1$ is a zero of $n_{\mathcal E}(s)^{-1}$ with $\Lambda(E,\lambda) \not=0$,
and $m_\lambda =2$
\item $\lambda =1$ and $-2-2J \le m_\lambda-2-2J \le$ the multiplicity of zero of $\Lambda(E,s)^2$ at $s=\lambda$,
where $J$ is the number of singular fibers of $\mathcal E$ (see \eqref{n_E}).
\end{enumerate}
Hence $w \ast h_{\mathcal E} \equiv 0$.
Because $w$ was arbitrary, we obtain ${\mathcal W}_\pi \subset {\mathcal T}(h_{\mathcal E})^{\bot}$.

The other implication ${\mathcal T}(h_{\mathcal E}) \subset {\mathcal W}_\pi^{\perp}$ is proved by a similar way.
The following fact is useful for this direction (see \cite[section 13]{GoJa} and \cite[section 2.5]{So});
there exists finitely many $(\phi_\alpha,f_{\pi,\alpha}) \in S(\pi)$ such that
\[
\sum_\alpha \int_{0}^{\infty} {\frak E}(\phi_\alpha,f_{\pi,\alpha})(x) x^{s} \frac{dx}{x} = L(\pi,s-1/2).
\]

The final assertion for ${\mathcal T}(h_{\mathcal E})^\perp = {\mathcal W}_\pi$ is obvious from \eqref{401} and (1) $\sim$ (5).
For ${\mathcal T}(h_{\mathcal E}) = {\mathcal W}_\pi^{\perp}$ we note that
${\mathcal W}_\pi^{\perp}$ consists of $f_{\lambda,k}$ such that
$\lambda$ is a zero of $n_{\mathcal E}(s)^{-2}\Lambda(E,s)s^4(s-2)^4(s-1)^2$ and $k \le$ the multiplicity of $\lambda$ (\cite{So}).
If $f_{\lambda,k} \in {\mathcal T}(h_{\mathcal E})$
then $\lambda$ is a pole of order $\ge k$ of ${\mathsf MC}(h_{\mathcal E})$
by the general theory of mean-periodic function (e.g. \cite[Theorem in lecture 4]{Kah}).
Hence the cancelation can not occur when ${\mathcal T}(h_{\mathcal E})= {\mathcal W}^\perp$. \hfill $\Box$


\subsection{Proof of Theorem \ref{thm_02}}
This is proved similarly to the proof of Theorem \ref{thm_01}.
For ${\cal T}(h_{\mathcal E} \circ {\ss}) \subset {\cal W}_\pi^\perp$,
it is sufficient to prove that $(h_{\mathcal E} \circ {\ss}) \ast w =0 $ for any $w \in {\mathcal W}_\pi$.
By the expansion of $h_{\mathcal E}(x)$ in Theorem \ref{prop_01},
$h_{\mathcal E} \circ {\ss}$ is a series consisting of $f_{\lambda,k} \circ {\ss}$.
We have
\[
\aligned
(f_{\lambda,k} \circ {\ss}) \ast w(y)
& = \sum_{j=0}^k (-1)^j \binom kj|y|^{-\lambda}(\log |y|)^{k-j} \int_{G_{\Q} \backslash G_{\A}}
w(x) |x|^\lambda (\log |x|)^{j} dx.
\endaligned
\]
Here
\[
\int_{G_{\Q} \backslash G_{\A}}
w(x) |x|^\lambda (\log |x|)^{j} dx = \frac{d^j}{d\lambda^j}\int_{G_{\Q} \backslash G_{\A}}
w(x) |x|^\lambda dx,
\]
and
\[
\aligned
\int_{G_{\Q} \backslash G_{\A}}
w(x) |x|^\lambda dx
& =
\int_{\R_+^\times} w_0(x) x^\lambda \frac{dx}{x} \\
& \times \int_{G_{\Q} \backslash G_{\A}} {\frak E}(\phi_1)(x)\varphi_\pi(x) |x|^\lambda dx
\int_{G_{\Q} \backslash G_{\A}} {\frak E}(\phi_2)(x)\varphi_\pi(x) |x|^\lambda dx,
\endaligned
\]
since $|~|^\lambda$ is a multiplicative (quasi) character.
By Lemma 3.5 of ~\cite{De},
\[
\int_{G_{\Q} \backslash G_{\A}} {\frak E}(\phi)(x)\varphi_\pi(x) |x|^\lambda dx
= L(\pi,s-1/2) F_{\phi,\varphi_\pi}(s) = L(E,s)F_{\phi,\varphi_\pi}(s),
\]
where $F_{\phi,\varphi_\pi}(s)$ is an entire function.
Therefore $(f_{\lambda,k} \circ {\ss}) \ast w(y)=0$ for each $\lambda$, $1 \le k \le m_\lambda$
appearing in the expansion of $h_{\mathcal E}$,
since $\lambda$ is a zero of $L(E,s)$ or a zero of $\int_{\R_+^\times} w_0(x) x^\lambda \frac{dx}{x}$.
Hence $(h_{\mathcal E} \circ {\ss}) \ast w =0 $ for any $w \in {\mathcal W}_\pi$.
\hfill $\Box$
%
%
%
\bibliographystyle{amsplain}
\bibliography{biblio}
\bigskip

\noindent
Masatoshi Suzuki \\
Department of Mathematics \\
Rikkyo University \\
Nishi-Ikebukuro, Toshima-ku \\
Tokyo 171-8501, Japan \\
\texttt{suzuki@@rkmath.rikkyo.ac.jp}

\end{document}